\documentclass[a4paper,12pt,DIV12]{scrartcl}
\usepackage{amsmath,amsthm,amssymb,enumerate}

\newtheorem{theorem}{Theorem}[section]
\newtheorem{proposition}[theorem]{Proposition}
\newtheorem{lemma}[theorem]{Lemma}

\newtheorem{definition}[theorem]{Definition}

\newenvironment{acknowledgements}{\bigskip\noindent\textbf{Acknowledgement}
\hskip 1mm}{\smallskip}

\DeclareMathOperator{\Aut}{Aut}

\DeclareMathOperator{\mlt}{Mlt} 
\DeclareMathOperator{\lmlt}{LMlt} \DeclareMathOperator{\rmlt}{RMlt}

\title{A class of simple proper Bol loops\thanks{This paper was written during
the author's Marie Curie Fellowship MEIF-CT-2006-041105.}}
\author{G\'abor P. Nagy}

\begin{document}

\maketitle

\begin{abstract}
The existence of finite simple non-Moufang Bol loops was considered as one of
the main open problems in the theory of loops and quasigroups. In this paper, we
present a class of proper simple Bol loops. This class also contains finite and
new infinite simple proper Bol loops. 
\end{abstract}

\section{Preliminaries}

For a loop $Q$, we call the maps $L_a(x)=ax, R_a(x)=xa$ \emph{left} and
\emph{right translations,} respectively. These are permutations of $Q$,
generating the \emph{left} and \emph{right multiplication groups} $\lmlt(Q),
\rmlt(Q)$ of $Q$, respectively. The group closure $\mlt(Q)$ of $\lmlt(Q)$ and
$\rmlt(Q)$ is the \emph{full multiplication group} of $Q$. Just like for groups,
normal subloops are kernels of homomorphisms of loops. The loop $Q$ is
\emph{simple} if it has no proper normal subloop. The
\emph{commutator-associator subloop} $Q'$ is the smallest normal subloop of $Q$
such that $Q/Q'$ is an Abelian group. For basic introductory reference on loops
see \cite{Pflugfelder}. 

The loop $Q$ is a \emph{left (right) Bol loop} if the identity
$x(y(xz))=(x(yx))z$ ($((xy)z)y=x((yz)y)$) holds in $Q$. Loops which satisfy
both identities are called \emph{Moufang loops.} 

For any field $F$, L.~J.~Paige \cite{Paige} constructed a simple nonassociative
Moufang loop $M(F)$. Using the classification of finite simple groups,
M.~Liebeck \cite{Liebeck} showed that the only finite simple nonassociative
Moufang loops are $M(\mathbb{F}_q)$. The existence of finite simple non-Moufang
Bol loops was considered as the one of the main open problems in the theory of
loops and quasigroups, cf. \cite{Problem} and \cite[Question 4]{Aschb}.

In this paper, Bol loops are left Bol loops and with proper Bol loops we mean
left Bol loops which are not Moufang.

\section{Bol loops as sections in groups}

Let $Q$ be a Bol loop, $G$ the group generated by the left translations of $Q$,
$H$ the stabilizer of $1$ in $G$ and $S$ the set of  left translations. ($S$ is
called the \emph{left section} of the loop.) It is widely known that the triple
$(G,H,S)$ satisfies the following properties: 

\begin{enumerate} 
\item[(F1)] $1\in S$ and for all $s,t \in S$, we have $sts \in S$. 
\item[(F2)] $S$ is a system of left coset representative for all conjugates
$H^g$ of $H$ in $G$. 
\end{enumerate}

Conversely, let $(G^*,H^*,S^*)$ be a triple satisfying (F1) and (F2). For $s,t
\in S^*$, let $u$ be the unique element in $stH^*\cap S^*$. Then, the binary
operation $s\circ t=u$ turns $S^*$ into a left Bol loop. Following M. Aschbacher
\cite{Aschb}, we call the triples $(G^*,H^*,S^*)$ satisfying (F1) and (F2)
\emph{Bol loop folders}.

The map $\alpha:S\to G/H$, $\alpha(s)=sH$ is a bijection. Let us denote by
$\tilde g$ the permutation induced by $g\in G$ on the left cosets of $H$ in
$G$: $\tilde g(hH)=ghH$. Then $\alpha L_s=\tilde s \alpha$, where $L_s:S\to S$,
$L_s(t)=s\circ t$ is the left translation in the Bol loop $(S,\circ)$. This
means that we have a surjective homomorphism 
\[\varphi: \hat G=\langle S \rangle \to \lmlt(S,\circ), \hskip 5mm s \mapsto
L_s=\alpha^{-1} \tilde s \alpha.\]
The kernel of $\varphi$ is the largest normal subgroup of $\hat G$ contained in
$H \cap \hat G$. If the permutation representation of $G$ on the left cosets of
$H$ is faithful, that is, if $H$ contains no proper normal subgroup of $G$,
then $\varphi$ is a bijection. 

In the remaining of this section, we show some \emph{folklore} results which
connect the structure of the Bol loop $Q$ and its left multiplication group. 

\begin{lemma} \label{lm:l=l'}
Let $Q$ be a loop with left multiplication group $G=\lmlt(Q)$. Then $Q=Q'$ if
and only if $G'$ acts transitively on $Q$. 
\end{lemma}
\proof Let $A$ be an Abelian group and assume $\varphi:Q\to A$ is a surjective
homomorphism. Then the map $\tilde \varphi: L_x\mapsto \varphi(x)$ extends to a
surjective homomorphism $G \to A$. Indeed, if $L_{x_1}\cdots L_{x_n}=1$ then
$x_1(x_2(\cdots x_n))=1$ and 
\[\tilde\varphi(L_{x_1})\cdots \tilde\varphi(L_{x_n}) =
\varphi(x_1)\cdots\varphi(x_n)=\varphi(x_1(x_2(\cdots x_n)))=\varphi(1)=1.\]
Moreover, $N=\ker\tilde\varphi$ contains the stabilizer $G_1$ of the unit
element of $Q$ in $G$:
\[\tilde\varphi(L_{xy}^{-1}L_xL_y)=\varphi(xy)^{-1}\varphi(x)\varphi(y)=1.\]
Therefore, $N$ and $G'\leq N$ are not transitive. Conversely, if $G'$ is not
transitive then $N=G'G_1$ is a proper normal subgroup of $G$. Then the map
$\varphi:Q\to G/N$, $\varphi(x)=L_xN$ is a surjective homomorphism:
\[ \varphi(x)\varphi(y)=L_xNL_yN=L_{xy}L_{xy}^{-1}L_xL_yN=L_{xy}N=\varphi(xy).
\hskip 5mm\qed  \]

\begin{lemma} \label{lm:nonsimpl_cond}
Let $Q$ be a Bol loop and let $\sigma$ be an automorphism of $\lmlt(Q)$ such
that $L_x^\sigma=L_x^{-1}$ for all $x \in Q$. Let $K$ be proper normal subloop
of $Q$. Then, $K=N(1)$ for some $\sigma$-invariant normal subgroup $N$ of
$\lmlt(Q)$. In particular, $Q$ is simple if all $\sigma$-invariant normal
subgroups of $\lmlt(Q)$ act transitively on $Q$. 
\end{lemma}
\proof Put $G=\lmlt(Q)$ and define the subset
\[ M=\{g \in G \mid g(yK)=yK \mbox{ for all $y\in Q$ }\}. \]
of $G$. Clearly, $M\lhd G$ and $L_x \in M \cap M^\sigma$ for all $x\in
K$. For the $\sigma$-invariant normal subgroup $N=M\cap M^\sigma$, $K=N(1)$
holds. \qed

We notice that for any simple proper Bol loop $Q$, the left multiplication
group possesses an involutorial automorphism $\sigma$ with $L_x^\sigma=L_x^{-1}$
for all $x \in Q$.

The bijection $u:Q\to Q$ is a \emph{left pseudo-automorphism} of the loop $Q$
with \emph{companion element} $c\in Q$ if 
\[cu(xy)=(cu(x))u(y)\] 
holds for all $x,y \in Q$. Equivalently, $u(1)=1$ and $L_c uL_x=L_{cu(x)} u$ for
all $x\in Q$. 

Two loops $(Q,\cdot)$ and $(K,\circ)$ are \emph{isotopes} if bijections
$\alpha,\beta,\gamma: Q\to K$ exist such that $\alpha(x)\circ
\beta(y)=\gamma(x\cdot y)$ for all $x,y\in Q$. The loop $Q$ is a \emph{G-loop}
if it is isomorphic to all its isotopes. 

\begin{lemma} \label{lm:leftcomp}
Let $Q$ be a left Bol loop and let us denote by $S$ the set of left
translations of $Q$. The loop $Q$ is a G-loop if and only if for all
$c\in Q$ there is an permutation $u$ of $Q$ such that $u(1)=1$ and
$uSu^{-1}=L_c^{-1}S$ hold.
\end{lemma}
\proof Combining Theorem III.6.1 and IV.6.16 of \cite{Pflugfelder}, we see that
a left Bol loop is a G-loop if and only if every element of $Q$ occurs as a
companion of some left pseudo-automorphism. Let $c\in Q$ be given and let
$u:Q\to Q$ be a bijection such that $u(1)=1$ and for all $x\in Q$ holds
$uL_xu^{-1}=L_c^{-1}L_{x'}$, where $x'\in Q$ depends on $x,c$ and $u$. Then,
$x'=cu(x)$ and $u$ is a left pseudo-automorphism with companion $c$. \qed

We close this section with a lemma on the left multiplication groups of
nonproper Bol loops.

\begin{lemma} \label{lm:simpmouf_lmlt}
Let $Q$ be a simple Moufang loop. Then $\lmlt(Q)$ is a simple group. 
\end{lemma}
\proof If $Q$ is a simple group, then $\lmlt(Q_{\mathcal T})\cong
Q_{\mathcal T}$ is simple. The left and right Bol identities can be written in
the form 
\[ R_{xz}R_x^{-1}=L_x^{-1}R_zL_x, \hskip 5mm L_{xy}L_y^{-1}=R_y^{-1}L_xR_y. \]
This means that for Moufang loops, the left and right multiplication
groups are normal in the full multiplication group. Theorem 4.3 of
\cite{NagyValsecchi} says that for an arbitrary nonassociative simple Moufang
loop $Q$, the multiplication group is simple. Hence, $\lmlt(Q)=\rmlt(Q)=\mlt(Q)$
is a simple group. \qed

\section{Construction of left Bol loops using extact factorizations of groups}

\begin{definition}
The triple $(X,Y_0,Y_1)$ is called an \emph{exact factorization triple} if $X$
is a group, $Y_0, Y_1$ are subgroups of $X$ satisfying $Y_0\cap Y_1=1$ and
$Y_0Y_1=X$. The exact factorization triple $(X,Y_0,Y_1)$ is \emph{faithful} if
$Y_0,Y_1$ do not contain proper normal subgroups of $X$.
\end{definition}

If $Y_1$ does not contain any proper normal subgroup of $X$, then an equivalent
definition of exact factorization triples is that $Y_0$ is a regular subgroup in
the permutation represetation of $X$ on the cosets of $Y_1$. In the
mathematical literature, the group $X$ is also called the \emph{Zappa-Sz\'ep
product} of the subgroups $Y_0,Y_1$. 

\begin{proposition} \label{pr:constr}
Let $\mathcal T=(X,Y_0,Y_1)$ be a faithful exact factorization triple. Let us
define the triple $(G, H, S)$ by 
\[ G=X\times X, \hskip 5mm H=Y_0 \times Y_1 \leq G, \hskip 5mm 
S= \{(x,x^{-1}) \mid x \in X\}. \]
Then $(G, H, S)$ is a Bol loop folder. The
associated left Bol loop $(S,\circ)$ is a G-loop. 
\end{proposition}
\proof We first show that $S$ is a left transversal for all conjugate of $H$.
Let
$a,b \in X$ be arbitrary elements; we can write $a^{-1}=a_0a_1$ and $b=b_0b_1$
in a unique way with $a_0,b_0 \in Y_0$ and $a_1,b_1 \in Y_1$. We have
\begin{eqnarray*}
\exists x \in X: (x,x^{-1}) \in (a,b)H &\Longleftrightarrow& \exists y_0
\in Y_0, y_1 \in Y_1: ay_0=y_1^{-1}b^{-1} \\
&\Longleftrightarrow& \exists y_0 \in Y_0, y_1 \in Y_1: a_1^{-1}a_0^{-1}y_0 =
y_1^{-1}b_1^{-1}b_0^{-1} \\
&\Longleftrightarrow& \exists y_0 \in Y_0, y_1 \in Y_1: a_0^{-1}y_0 b_0 =
a_1y_1^{-1}b_1^{-1} \in Y_0 \cap Y_1.
\end{eqnarray*}
Since $Y_0\cap Y_1=1$, we obtain $y_0=a_0b_0^{-1}$, $y_1=b_1^{-1}a_1$ and the
unique element of $(a,b)H \cap S$ is $(a_1^{-1}b_0^{-1}, b_0a_1)$. This shows
that $S$ is a left transversal to $H$ in $G$. In order to prove the same fact
for the conjugates of $H$, let us take an arbitrary $g \in G$ and write $g=sh$;
we have $H^{g}=H^{s}$. For $a \in G$ let us define $t\in S$ as the unique
element of $S\cap s^{-1}as^{-1}H$. Then $sts$ is the unique element of $S\cap
aH^{s}$. This proves (F1) and (F2).

Let $b\in X$ and write $b=y_0y_1^{-1}$ with $y_0 \in
Y_0, y_1 \in Y_1$. Then
\[(y_0,y_1)(a,a^{-1})(y_0,y_1)^{-1} = (b,b^{-1})(y_1 a y_0^{-1}, y_0 a^{-1}
y_1^{-1}) \in (b,b^{-1}) S.\]
Since $b\in X$ is arbitrary and $(y_0,y_1) \in H$, Lemma \ref{lm:leftcomp}
implies that $(S,\circ)$ is a G-loop. \qed

\begin{definition}
Let $\mathcal T=(X,Y_0,Y_1)$ be a faithful exact factorization triple and let us
define $G,H,S$ as in Proposition \ref{pr:constr}. The Bol loop corresponding to
the Bol loop folder $(G, H, S)$ will be denoted by $Q_{\mathcal T}$. 
\end{definition}

\begin{lemma} \label{lm:lmlt}
Let $\mathcal T=(X,Y_0,Y_1)$ be a faithful exact factorization triple and let us
define $G,H,S$ as in Proposition \ref{pr:constr}. Then $\lmlt(Q_{\mathcal T})$
is isomorphic to $\hat G = \langle S \rangle \cong X.X' \lhd G$. 
\end{lemma}
\proof We claim that $H$ contains no normal subgroup of $G$. Indeed, the
projections of $H$ to the direct factors of $G=X \times X$ are $Y_0,Y_1$ which
contain no normal subgroup of $X$. Thus a normal subgroup of $H$ must have
trivial projections, hence it must be trivial. Therefore, the permutation action
of $G$ on the left cosets of $H$ is faithful and we can consider $G$ as a
permutation group. Moreover, by the definition of $Q=Q_{\mathcal T}$, the left
translations are precisely the permutations induced by the elements of $S$. This
proves the lemma. \qed

\begin{proposition} \label{pr:nonsolv_cond}
Let $\mathcal T=(X,Y_0,Y_1)$ be a faithful exact factorization triple such that
$X'Y_0=X''Y_1=X$ and let us define $G,H,S$ as in Proposition \ref{pr:constr}.
Then 
\begin{enumerate}
\item[(i)] $\hat G'H=G$ with $\hat G=\langle S \rangle$. 
\item[(ii)] $Q'_{\mathcal T}=Q_{\mathcal T}$. In particular, $Q_{\mathcal T}$ is
not solvable.
\end{enumerate}
\end{proposition}
\proof Since $X'\times X' \leq \hat G$, we have $X'' \times X'' \leq \hat G'$.
By $1\times Y_1 \leq H$ and $X''Y_1=X$, we obtain $1 \times X \leq \hat G'H$.
Clearly, for any $a,b \in X$, $([a,b],[a^{-1},b^{-1}]) \in \hat G'$, thus
$([a,b],1) \in \hat G'H$. This implies $X'\times 1 \leq \hat G'H$, hence
$X\times 1 \leq \hat G'H$ by $Y_0\times 1\leq H$ and $X'Y_0=X$. This proves
(i). Lemma \ref{lm:lmlt} says that $\hat G$ is isomorphic to the left
multiplication group of $Q$ and the commutator subgroup of $\hat G$ acts
transitively on the left cosets of $H$ in $G$. Therefore, (ii) follows from
Lemma \ref{lm:l=l'}. \qed

We call the group $X$ \emph{almost simple} if $T\leq X \leq \Aut(T)$ for some
nonabelian simple group $T$. The group $T$ is the \emph{socle} of $X$.

\begin{theorem} \label{th:simple_cond1}
Let $X$ be an almost simple group with socle $T$. Let $\mathcal T=(X,Y_0,Y_1)$
be a faithful exact factorization triple and assume $X=TY_0=TY_1$. Then
$Q_{\mathcal T}$ is a simple proper left Bol loop. 
\end{theorem}
\proof Let $\sigma$ be the automorphism of $G$ mapping $(a,b)\mapsto (b,a)$.
Since
$S^\sigma=S$, we have $\hat G^\sigma=\hat G$ and $\lmlt(Q_{\mathcal T})$
has an automorphism which inverts the left translations. Clearly, $T\times T
\leq X' \times X' \leq \hat G$ and every $\sigma$-invariant normal subgroup of
$\hat G$ contains $T\times T$. However, $T\times T$ is transitive by
assumption, thus, $Q_{\mathcal T}$ is simple by Lemma \ref{lm:nonsimpl_cond}.
Moreover, $Q$ is proper Bol by Lemma \ref{lm:simpmouf_lmlt}. \qed

\section{Some classes of simple proper Bol loops}

In this section we present some finite and infinite simple proper Bol loops by
applying the construction of Proposition \ref{pr:constr}.

\medskip\noindent\textbf{Example I:} Put $X=PSL(n,2)$, let $Y_0$ be a Singer
cycle and $Y_1$ be the stabilizer of a projective point. Then $Q_{(X,Y_0,Y_1)}$
is a finite simple proper Bol loop by Theorem \ref{th:simple_cond1}. We notice
that many other finite simple groups have exact factorizations. The
factorizations of finite groups are intensively studied, cf. \cite{LiePrChS},
\cite{Giudici} and the references therein.

\medskip\noindent\textbf{Example II:} Let $n$ be an even integer and put
$X=S_n$, $Y_0=\langle (1,2,\ldots,n) \rangle$ and $Y_1=S_{n-1}$ with $n\geq 4$.
Define the loop $Q_n=Q_{(X,Y_0,Y_1)}$. If $n\geq 6$ then $Q_n$ is simple by
Theorem \ref{th:simple_cond1}. If $n=4$ then by Proposition
\ref{pr:nonsolv_cond} $Q_4$ is a nonsolvable Bol loop of order $24$. It is known
that all Bol loops of order at most $12$ are solvable, thus, $Q_4$ is simple. We
emphasize the fact that the left multiplication group of $Q_4$ is a solvable
group of order $288$. The computer result \cite{Moorh} of G.~E.~Moorhouse shows
that all Bol loops of order less than $24$ are solvable, hence $Q_4$ is a
simple Bol loop of least possible order.

\medskip\noindent\textbf{Example III:} Put $X=PSL_2(\mathbb R)$ and define
the subgroups 
\[Y_0=\left\{\pm \left ( \begin{array}{cc} \cos t &\sin t \\ -\sin t&\cos t
\end{array} \right ) \mid t \in \mathbb R \right\}, \hskip 5mm 
Y_1=\left \{\pm \left ( \begin{array}{cc} a &b \\ 0&a^{-1} \end{array} \right )
\mid  a \in \mathbb R \setminus \{0\}, b \in \mathbb R \right \}\]
of $X$. By Theorem \ref{th:simple_cond1}, $Q_{(X,Y_0,Y_1)}$ is a proper 
simple proper Bol loop which is isomorphic to all its isotopes. In particular,
$Q_{(X,Y_0,Y_1)}$ is not isotopic to a Bruck loop. Moreover, the left
translation group is $PSL_2(\mathbb R) \times PSL_2(\mathbb R)$. 

In \cite{Figula}, the author classifies all differentiable Bol loops having a
semi-simple left multiplications group of dimension at most $9$. Our
construction shows that the classification cannot be complete. (The author seems
not to consider the case when the group $G$ topologically generated by the left
translations is a proper direct product of simple Lie groups $G_1,G_2$ and the
stabilizer of $1\in Q$ in $G$ is a direct product $H=H_1\times H_2$ with
$1\neq H_i\leq G_i$, $i=1,2$.)

\medskip\noindent\textbf{Example IV:} Let $\Sigma$ be the set of non-zero
squares in $\mathbb F_{27}$, $|\Sigma|=13$. Let $X$ be the set of
transformations
\[X=\{ f(z)=az^\tau+b \mid a \in \Sigma, b \in \mathbb F_{27}, \tau \in
\Aut(\mathbb F_{27}) \}\]
of $\mathbb F_{27}$. Then $X$ has order $1053=3^4\cdot 13$ and it acts
primitively on $\mathbb F_{27}$. Moreover, 
\[X'=\{ f(z)=az+b \mid a\in \Sigma, b\in \mathbb F_{27}\}, \hskip 5mm
X''=\{ f(z)=z+b \mid b\in \mathbb F_{27}\}. \]

We define $Y_1$ as the stabilizer of $0$ in $X$. Since $X''$ acts regularly, we
have $X''Y_1=X$. Let $U$ be the $3$-Sylow subgroup of $X$, $|U|=81$. Clearly,
$U/X''\cong C_3$, thus $U'\cap U_0=1$ where $U_0$ is the stabilizer of $0$ in
$U$. Therefore, $U$ has a subgroup $Y_0$ of order $27$ such that $Y_0\neq X''$
and $Y_0 \cap U_0=1$; in other words, $Y_0$ acts regularly on $\mathbb F_{27}$.
Since $X''$ is the unique $3$-Sylow subgroup of $X'$, $Y_0$ cannot be contained
in $X'$. This implies $X'Y_0=X$ because $X'$ has index $3$ in $X$. 

We define now the Bol loop $Q=Q_{(X,Y_0,Y_1)}$. Let $K$ be a maximal proper
normal subloop of $Q$, that is, $Q/K$ be a simple loop. If $Q/K$ were
associative then by the Odd Order Theorem, it would be a cyclic group and we
had a surjective homomorphism from $Q$ to an Abelian group. By Proposition
\ref{pr:nonsolv_cond} this is not possible. Hence, $Q/K$ is a proper simple Bol
loop of odd order. (It can be shown by computer that $Q$ itself is a simple Bol
loop.) This last construction shows that the Odd Order Theorem does not hold for
finite Bol loops. (Cf. \cite{Twisted}.)

\begin{acknowledgements}
I would like to thank Peter M\"uller (Uni. W\"urzburg) for his help in finding
the group $X$ in Example IV. I also thank Petr Vojt\v{e}chovsk\'y and
Michael Kinyon (Uni. Denver) for many stimulating conversations and helpful
comments.
\end{acknowledgements}

\end{document}